\documentclass[12pt]{article}
\usepackage{amssymb}

\author{I. G. Korepanov\\
\normalsize Southern Ural State University\\[-0.5ex]
\normalsize 76 Lenin avenue\\[-0.5ex]
\normalsize 454080 Chelyabinsk, Russia\\[-0.5ex]
\normalsize E-mail: kig@susu.ac.ru}
\date{}
\title{Euclidean 4-simplices and invariants of four-dimensional
manifolds: \goodbreak III.~Moves~$1\leftrightarrow 5$ and related structures}

\def\be{\begin{equation}}
\def\ee{\end{equation}}

\def\pa#1#2{{\partial#1\over\partial#2}}

\def\Im{\mathop{\rm Im}\nolimits}

\def\minor{\mathop{\rm minor}\nolimits}

\def\emline#1#2#3#4#5#6{%
       \put(#1,#2){\special{em:moveto}}%
       \put(#4,#5){\special{em:lineto}}}
\def\newpic#1{}

\begin{document}
\maketitle

\begin{abstract}
We conclude the construction of the algebraic complex, consisting of spaces
of differentials of Euclidean metric values, for four-dimensional
piecewise-linear manifolds. Assuming that the complex is acyclic, we
investigate how its torsion changes under rebuildings of the manifold
triangulation. First, we write out formulas for moves $3\to 3$
and~$2\leftrightarrow 4$ based on the results of our two previous works,
and then we study in detail moves~$1\leftrightarrow 5$. On this basis,
we obtain the formula for a four-dimensional manifold invariant.

As an example, we present a detailed calculation of our invariant
for sphere~$S^4$; in particular, the complex turns out, indeed, to be acyclic.
\end{abstract}


\section{Introduction}
\label{sec IIIvvedenie}

This work is third in the series of papers started with papers \cite{I}
and~\cite{II}. We also use Roman numerals I and~II, respectively, for
references to those works; in particular, equation~(II.1) is the
equation~(1) from paper~\cite{II}, or ``from paper~II''.

Our goal is to construct and investigate a new type of acyclic complexes,
wherefrom we should be able to extract invariants of four-dimensional
piecewise-linear manifolds. Note, however, that at the time when paper~I
was being written, it was not yet clear that acyclic complexes were exactly
those structures that stood behind the considered algebraic formulas,
as well as behind the invariants of {\em three-dimensional\/} manifolds from
papers \cite{3dim1} and~\cite{3dim2}. These complexes were written out
in paper~II --- the whole complex for the three-dimensional case and a part
of it --- for the four-dimensional case. Concerning the four-dimensional
case, it turned out that, in order to construct a value invariant
under Pachner moves (i.e., rebuildings of the manifold triangulation) of
type $3\to 3$ {\em only\/}, it is enough to consider only a small central part of the
complex (two vector spaces and one linear mapping between them), whereas
the addition of rebuildings $2\leftrightarrow 4$ brings into
consideration one more space (namely, the space of edge deviations) and
one more mapping.

Paper~I was devoted to moves~$3\to 3$, and paper~II ---
to moves~$2\leftrightarrow 4$. Therefore, the complex {\em in full\/} was
not necessary in those papers, and its full form was only announced
in paper~II (formulas~(II.6) and~(II.7)).

The contents of the present ``paper~III'' by sections is as follows.
In Section~\ref{sec IIIp2} we revisit, once again, the dimensionality~3
and write out the acyclic complex from Section~2 of paper~II in a slightly
different (longer) form, which clarifies the structure of one of the
involved spaces as a {\em factor\/} space. This simple improvement
may simplify considerably our constructions when we pass to higher
dimensions.

In Section~\ref{sec IIIp3} we present the space of vertex deviations,
as we promised in paper~II, and its mapping into the space of edge
deviations. This enables us to write out in Section~\ref{sec IIIp4}
our complex for the four-dimensional case in full; we also ``make it longer''
with respect to paper~II in analogy with the three-dimensional case.
We still have to prove, though, that our sequence of spaces and mappings
is really a complex, and we do that in Section~\ref{sec IIIp5}.

We would like our constructed complex to be acyclic (i.e., an exact
sequence) which would enable us to get manifold invariants out of its
torsion. We do not give the full proof of the acyclicity in the general case,
limiting ourselves to some remarks in the end of Section~\ref{sec IIIp5}
(see also a concrete example in Section~\ref{sec IIIp9}). Then we show
that, {\em assuming the exactness}, the torsion of our complex does not change
under Pachner moves (and the exactness itself is conserved under those moves).
To this end, we study in Section~\ref{sec IIIp6} how the complex changes
under moves $3\to 3$ and~$2\leftrightarrow 4$, and in Section~\ref{sec IIIp7}
 --- under moves~$1\leftrightarrow 5$ (one can say that it was the need
to study the moves~$1\leftrightarrow 5$ that required the detailed
consideration of the {\em whole\/} algebraic complex).
In Section~\ref{sec IIIp8} we present the final formula giving the invariant
of a four-dimensional piecewise-linear manifold (i.e., the invariant of
{\em all\/} Pachner moves) in terms of the torsion of the complex and other
Euclidean geometric values. In Section~\ref{sec IIIp9} we show how our
machinery works for the sphere~$S^4$ (and check explicitly that, at least
for the sphere, our complex is indeed acyclic).

In the concluding Section~\ref{sec IIIobsuzhdenie} we discuss the obtained
results.

\section{Revisiting three dimensions once again: a longer form of the acyclic
complex}
\label{sec IIIp2}

We recall that the acyclic complex from Section~2 of paper~II,
corresponding to a triangulation of a three-dimensional manifold~$M$,
had the following form:
$$
0
\leftarrow (\cdots)
\stackrel{B^{\rm T}}{\leftarrow} (d\omega) \stackrel{A}{\leftarrow} (dl)
\stackrel{B}{\leftarrow} (dx
\hbox{ and } dg) \leftarrow 0.
\eqno\hbox{(II.1)}
$$

Here $(d\omega)$ is the space of column vectors made of infinitesimal
deficit angles at all edges of the complex, $(dl)$ is the space of
columns of differentials of edge lengths, $(dx)$ is the space of columns of
differentials of Euclidean coordinates of the edges of the complex
(to be exact, of their lift-ups onto the universal cover), $(dg)$ is the space
of columns of differentials of continuous paprameters on which the
representation $f\colon \; \pi_1(M)\to E_3$ may depend. In particular,
$(dg)$ vanishes for manifolds with a finite fundamental group.

We thus continue to use the notations in the style of paper~II for
linear spaces entering in algebraic complexes, hoping that the
convenience of such notations pays for their certain looseness. Of course,
 from the formal standpoint, for instance, space~$(dl)$ is the {\em
tangent space\/} to the manifold consisting of all sets of positive
numbers --- ``lengths'' --- put in correspondence to the edges of our
triangulation; $(dx\hbox{ and } dg)$ (the space of columns of all~$dx$ and
all~$dg$) is the {\em direct sum\/}~$(dx)\oplus(dg)$.

Recall that $(dx)$ in sequence~(II.1) is the space of differentials of
coordinates taken up to those infinitesimal motions of the Euclidean
space~$\mathbb R^3$ which are compatible with the given representation~$f$,
i.e.\ motions commuting with its image $\Im f=f\bigl( \pi_1(M)\bigr)$. Such
motions form a subalgebra in the Lie algebra $\mathfrak e_3$, which we
denote~$\mathfrak a$. Column~$(dx)$ is described explicitly for various
cases in Section~2 of paper~II. In the three-dimensional case, such
description can be made easily, but difficultuies may increase when we pass to
higher dimensions. That is the reason for rewriting sequence~(II.1)
in the following, more elegant form.

We now permit ourselves to {\em change notations\/} and understand below
by $(dx)$ the space of columns of differentials of Euclidean coordinates of
the vertices in the complex, with no further conditions, that is,
columns of values $(dx_1,dy_1,dz_1,\ldots,\allowbreak dx_N,dy_N,dz_N)$.
Then our {\em former\/} space~$(dx)$ is written as the
factor $(dx)/\mathfrak a$. This suggests an idea of adding one more term
to the right of $(dx\hbox{ and }dg)$ in sequence~(II.1) (and a symmetric
term in the left-hand part of the sequence --- see below), understanding
now $(dx)$ in the new sense:
\be
0 \leftarrow (\cdots)
\leftarrow (\cdots)
\stackrel{B^{\rm T}}{\leftarrow} (d\omega) \stackrel{A}{\leftarrow} (dl)
\stackrel{B}{\leftarrow} (dx
\hbox{ and } dg) \leftarrow \mathfrak a \leftarrow 0.
\label{3 dlinnaya}
\ee
Recall that the whole sequence (this applies to both (II.1) and
(\ref{3 dlinnaya})) is symmetric in the following sense. Each term in it
is considered as a vector space with a fixed basis; mappings between them are
identified with matrices; and any two matrices at equal distances from the
left-hand and right-hand ends of the sequence must be obtained from one
another by matrix transposing. In particular, matrix~$A$ that gives the
mapping $(dl)\to (d\omega)$ is symmetric.

As for the notation ``$(\cdots)$'', we are using it for {\em different\/}
linear spaces whose specific geometrical sense we are not going to
investigate (at this moment; still, we know the matrices of mappings between
such spaces, for instance, from the symmetry described in the previous
paragraph).

The torsion of an acyclic complex is the product of certain minors in the
matrices of its linear mappings, taken in alternating degrees~$\pm 1$. If we
adopt a convention that the $+$~sign corresponds to the first nontrivial
mapping (coming after the zero injection, i.e., for example, in
sequence~(II.1) that is the second arrow from the right), then the torsion
will change to its inverse when we pass from (II.1) to~(\ref{3 dlinnaya}).
To avoid this, we can {\em agree\/} to take the $+$~sign for mapping~$B$ in
sequence~(\ref{3 dlinnaya}), as before. Then {\em the torsion of
complex~(\ref{3 dlinnaya}) coincides with the torsion of complex~(II.1)\/}
if we choose a basis in Lie algebra~$\mathfrak e_3$ in a natural way
(namely, three infinitesimal translations along mutually orthogonal axes and
three rotations around these very axes).

The full proof of this statement must involve the analysis of all the particular cases from
Section~2 of paper~II. Here, we will limit ourselves to the case of a lens
space $L(p,q)$, with a nontrivial homomorphism
$f\colon \; \pi_1\bigl(L(p,q)\bigr) \to E_3$ whose image consists of rotations
around the $z$~axis. To compose the minor of the matrix of mapping
$\mathfrak a \to (dx)$, we choose two basis elements in space~$(dx)$,
namely, $dx_1$ and $dz_1$, where subscript~$1$ corresponds to some vertex
in the triangulation which we agree to call the ``first'' one. Now a simple
straightforward calculation shows that, indeed, the torsions of complexes
(\ref{3 dlinnaya}) and (II.1) coincide, if we choose the basis of the
``old''~$(dx)$ in the latter complex according to formula~(II.2)
(it makes sense to remind here that we consider the torsion {\em to within
its sign\/}).

\section{Vertex deviations and their mapping into edge deviations}
\label{sec IIIp3}

We have already used Euclidean coordinates in our constructions a few times,
namely, the coordinates of a vertex in the complex were necessary to define
the mappings $(dx)\to (dl)$ in the three-dimensional and four-dimensional
cases, and the coordinates of the vector of edge deviation --- for the mapping
$(d\vec v)\to (dS)$. We can remark that we did not need any connection between
the coordinate systems for {\em different\/} vertices and/or edges while
doing these constructions. Any individual coordinate system could be chosen
arbitrarily, for instance, by fixing some angles between coordinate axes
and adjacent edges.

We cannot work with no such coordinate systems at all if we want to
{\em fix the bases\/} in all vector spaces entering in a complex; on the
other hand, the torsion of our complexes, as one can check, does not
depend on the choice of those systems.

In this Section we, first, define the {\em vertex deviation\/} as a tensor
value which needs, for its components to be fixed, a coordinate system
corresponding to the {\em vertex}. Next, it will be convenient for us to
define how it generates edge deviations for the edges abutting on the given
vertex, using the {\em same\/} (i.e.\ corresponding to the vertex) coordinate
system for those edges. We {\em imply} that the components of each edge's
deviation are then transformed into {\em its own\/} coordinate system by
using a proper orthogonal transformation.

So, we call vertex deviation a bivector (antisymmetric tensor)
$d\sigma_{\alpha \beta }$, where $\alpha, \beta=1,\ldots, 4$. Let $AB$ be
one of the edges abutting on vertex~$A$. This edge can be also considered as a
four-dimensional vector $\overrightarrow{AB}$, whose coordinates we denote
as~$l_{\alpha}$. By definition, a deviation of vertex~$A$ equal to
$d\sigma_{\alpha \beta }$ generates the deviation of edge $d\vec v_{AB}$
with components
\be
\left( d\vec v_{AB} \right)_{\beta} = \frac{1}{L_{AB}}\, \sum_{\alpha}
l_{\alpha}\, d\sigma_{\alpha \beta },
\label{sigma->v}
\ee
where $L_{AB}=\sum_{\alpha} l_{\alpha}^2$ is the squared length of edge~$AB$.

In general, $d\vec v_{AB}$ is the sum of expression (\ref{sigma->v}) and
a similar expression involving the deviation of vertex~$B$ (where, of course,
vector $\overrightarrow{BA}$ must be used instead of~$\overrightarrow{AB}$).
The reasonableness of definition~(\ref{sigma->v}) will be clear
in Section~\ref{sec IIIp5}, where we will prove, in particular, that
the edge deviations of type~(\ref{sigma->v}) generate zero
differentials~$(dS)$ of two-dimensional face areas. This will be part of
the statement that ``composition of two neighboring arrows (i.e.\ linear
mappings) in the sequence is zero'', which will justify the name ``complex''
for that sequence. Let us now pass on to the construction of our sequence
in its full form.

\section{The full sequence of spaces and mappings in the
four-dimensional case}
\label{sec IIIp4}

We will write out two ``conjugate'' sequences, as we did in Section 3 of
paper~II. Compared with formulas (II.6) and~(II.7), we will increase them in
length both on the right and on the left, in analogy with the three-dimensional case
(Section~\ref{sec IIIp2} of the present work), because both the leftmost
nonzero space in formula~(II.7) and the rightmost one in~(II.6) can be
represented naturally as factor spaces.

As for the first of mentioned spaces, denoted as $(dx\hbox{ and }{dg})$
in formula~(II.7), here our elongation goes, in principle, the same way as
in three dimensions. Thus we pass at once to the second one, denoted
$(d\sigma)$ in formula~(II.6), i.e.\ to the space of vertex deviations
introduced in the previous Section. It turns out that there exists
an easily defined space of ``trivial'' deviations generating zero edge
deviations~$d\vec v_a$, and it finds its natural place in our complex.

Consider formula (\ref{sigma->v}). We recall that, in it, the tensor
value~$d\sigma_{\alpha \beta }$ pertains to point~$A$, and $l_{\alpha}$ are
the components of vector~$\overrightarrow {AB}$. If we take into account
$d\sigma_{\alpha \beta }$ in point~$B$ as well, we get
\be
\left( d\vec v_{AB} \right)_{\beta} = \frac{1}{L_{AB}}\, \sum_{\alpha}
l_{\alpha} \left( (d\sigma_A)_{\alpha \beta } - (d\sigma_B)_{\alpha \beta }
\right).
\label{4.1}
\ee
One can see from here that $d\vec v_{AB}=0$ if the difference
$(d\sigma_A)_{\alpha \beta } - (d\sigma_B)_{\alpha \beta }$
has the form
$\displaystyle \sum_{\gamma,\delta} \epsilon_{\alpha \beta \gamma \delta}\,
l_{\gamma}\, ds_{\delta}$,
where $\epsilon_{\alpha \beta \gamma \delta}$ is the totally antisymmetric
tensor, $\epsilon_{1234}=1$, and $d\vec s$ is any infinitesimal vector.

Passing on from considering one edge $AB$ to considering all edges in the
complex, we take, first, the case of trivial representation
$f\colon \; \pi_1 (M) \to E_4$ of the fundamental group of our manifold~$M$
into the group of motions of the four-dimensional Euclidean space,
i.e.\ the case $\Im f =\{ e \}$. This means that all inverse images of any
given vertex in the complex get into one and the same point of space~$\mathbb R^4$
(see Section~3 of paper~I), that is, simply speaking, to each vertex~$A$
its radius vector~$\vec r_A$ corresponds unambiguously with components
$(r_A)_{\alpha}$ in some Cartesian coordinate system common for the whole
complex.

Choose some infinitesimal antisymmetric tensor $d\tau_{\alpha\beta}$ and
vector $ds_{\alpha}$, and set for each vertex~$A$
\be
(d\sigma_A)_{\alpha\beta} = d\tau_{\alpha\beta} +
\epsilon_{\alpha \beta \gamma \delta} (r_A)_{\gamma} \, ds_{\delta}.
\label{trivsigma}
\ee
It follows from the foregoing that such set of vertex deviations yields zero
deviations for all edges. Note that such columns $d\sigma$ make up a
ten-dimensional linear space.

In case representation $f$ is not trivial, any vertex has more than one
inverse images in the universal covering, and these inverse images are placed
in different points of space~$\mathbb R^4$. The requirement that
formula~(\ref{trivsigma}) must give equal results for all such
inverse images (after transforming them into the coordinate system corresponding to
the given vertex --- see Section~\ref{sec IIIp3}), leads to linear
restrictions on admissible $d\tau_{\alpha\beta}$ and~$ds_{\beta}$.
In the present paper, we do not write out explicitly these restrictions:
it will suffice for us to begin with studying just simply connected manifolds,
for which $\pi_1(M)=\{e\}$. Nevertheless, we introduce notation $(d\sigma)_0$
for the subspace of columns of those vertex deviations which are correctly
determined by formula~(\ref{trivsigma}), for {\em any\/} manifold~$M$.

Now we are ready to rewrite sequences (II.6) as~(II.7) in the renovated
(longer) form. Despite the fact that we are still using notation
``$(\cdots)$'' for vector spaces with whose geometric sense we are not
concerned now, we have {\em given the definitions\/} to all {\em matrices\/}
of mappings denoted by arrows. So, here are our mutually conjugate sequences:
\begin{eqnarray}
0\leftarrow (\cdots) \leftarrow (\cdots) \leftarrow (d\Omega_a) \stackrel{(\partial
\Omega_a / \partial S_i)}{\longleftarrow} (dS_i) \leftarrow (d\vec v_a)
\leftarrow (d\sigma) \leftarrow (d\sigma)_0 \leftarrow 0\,,
\label{posled} \\
0 \rightarrow \mathfrak a \rightarrow (dx \hbox{ and } dg) \rightarrow (dL_a) \stackrel{(\partial
\omega_i / \partial L_a)}{\longrightarrow} (d\omega_i) \rightarrow
(\cdots) \rightarrow (\cdots) \rightarrow (\cdots) \rightarrow 0\,.
\label{posledT}
\end{eqnarray}
Here, of course, $\mathfrak a$ is a subalgebra of Lie algebra $\mathfrak e_4$
of motions of Euclidean space~$\mathbb R^4$. Namely, $\mathfrak a$ consists of
those motions commuting with the image of group $\pi_1(M)$ in~$E_4$, in full
analogy with the three-dimensional case from Section~\ref{sec IIIp2}.

\section{The sequence is a complex}
\label{sec IIIp5}

Now we will show that the composition of any two successive arrows in
sequence~(\ref{posled}) or, equivalently, (\ref{posledT}) equals zero.
Thus, we will justify the name ``complex'' for each of these sequences.

We start with the two arrows of sequence~(\ref{posled}) adjacent to the
term~$(d\sigma)$. It follows directly from formulas (\ref{4.1})
and~(\ref{trivsigma}) of the previous Section that their composition is zero.

Moving to the left, we consider two arrows around the term~$(d\vec v_a)$.
Consider a triangle~$ABC$ --- one of the two-dimensional faces of our
simplicial complex --- and check that the deviations of its edges, if they
are given by formulas of type~(\ref{4.1}), lead to the zero area
differential~$dS_{ABC}$.

It is enough to consider the case where only vertex~$A$ has a nonzero
deviation~$d\sigma$. According to formula~(\ref{sigma->v}), the lengths
of vectors $d\vec v$ in Figure~\ref{IIIris1}
\begin{figure}
\begin{center}
\unitlength=0.75mm
\special{em:linewidth 0.5pt}
\linethickness{0.5pt}
\begin{picture}(66.00,66.00)
\emline{5.00}{5.00}{1}{65.00}{5.00}{2}
\emline{65.00}{5.00}{3}{25.00}{65.00}{4}
\emline{25.00}{65.00}{5}{5.00}{5.00}{6}
\put(45.00,35.00){\vector(-3,-2){6.00}}
\put(15.00,35.00){\vector(-3,1){7.80}}
\special{em:linewidth 0.2pt}
\linethickness{0.2pt}
\emline{7.20}{37.60}{7}{35.00}{5.00}{8}
\emline{35.00}{5.00}{9}{39.00}{31.00}{10}
\emline{39.00}{31.00}{11}{7.20}{37.60}{12}
\put(15.00,38.50){\makebox(0,0)[rb]{$d\vec v_{AB}$}}
\put(7.00,36.00){\makebox(0,0)[rt]{$C'$}}
\put(42.00,34.00){\makebox(0,0)[rb]{$d\vec v_{AC}$}}
\put(39.50,30.50){\makebox(0,0)[lt]{$B'$}}
\put(35.00,4.00){\makebox(0,0)[ct]{$A'$}}
\put(4.00,4.00){\makebox(0,0)[rt]{$B$}}
\put(66.00,4.00){\makebox(0,0)[lt]{$C$}}
\put(25.00,66.00){\makebox(0,0)[cb]{$A$}}
\end{picture}
\end{center}
\caption{}
\label{IIIris1}
\end{figure}
are inversely proportional to the trianle's sides to which they belong.
Thus, one can easily deduce that the area of triangle~$A'B'C'$ multiplied by four
{\em coincides\/} with the area of triangle~$ABC$ and hence $dS_{ABC}=0$
according to formula~(II.8).

Moving further to the left along sequence~(\ref{posled}), we must consider
two arrows adjacent to the term $(dS_i)$ --- but the vanishing of their
product has been already proven in Section~5 of paper~II.

To investigate the two remaining pairs of arrows, we switch to
sequence~(\ref{posledT}). In terms of that sequence, these are the pairs of
arrows around the terms $(dx \hbox{ and } dg)$ and~$(dL_a)$. Now it remains to
remark that the statements we need can be proved in full analogy to how it was
done for the three-dimensional case in Section~2 of paper~II and Section~2
of the present work. Moreover, there is even {\em exactness\/} in those
terms.

As for the exactness in other terms, nothing is known about it as yet in the
general case, save that it is clear from the construction that the exactness holds
in the left- and rightmost terms $\mathfrak a$ and~$(d\sigma)_0$ which appeared
when we ``made longer'' our complex. Still, we will see in the following
Sections \ref{sec IIIp6} and~\ref{sec IIIp7} that if the sequence is exact,
then this property is preserved under Pachner moves, that is, exactness does not
depend on a triangulation. Besides, the example of sphere~$S^4$ studied
below in Section~\ref{sec IIIp9} shows that at least for the sphere the
sequence {\em is\/} exact (i.e., is an acyclic complex).

\section{How the algebraic complex changes under moves $3\to 3$
and~$2\leftrightarrow 4$}
\label{sec IIIp6}

We already know from papers I and~II what happens under moves $3\to 3$
and~$2\leftrightarrow 4$. Our task in this Section is to reformulate these
results while holding strictly to the algebraic language of acyclic complexes.

\subsection{Moves $3\to 3$}
\label{subsec 6.1}

In paper~I, devoted to moves~$3\to 3$, we were considering matrix
$(\partial \Omega_a / \partial S_i)$ and its conjugate
$(\partial \omega_i / \partial L_a)$. From our current viewpoint, they form
the central part of sequence (\ref{posled}) and its conjugate~(\ref{posledT}).
Some of the results of paper~I (namely, Theorem~4) do not deal with moves
$3\to 3$ as such but only reproduce (by somewhat amateurish means) a part
of statement that the torsion of a complex does not depend on a specific
choice of minors through which it is expressed. It is Theorem~3 of paper~I
that deals with moves $3\to 3$ proper. It seems worthwhile to recall it here
once again (with only a very slight reformulation):

\medskip

{\it
Choose in matrix $(\partial \omega_i/ \partial L_a)$ a largest
square submatrix~$\cal B$ with nonzero determinant. Let
$\cal B$ contain a row corresponding to such face $i=ABC$ that
belongs to exactly three 4-simplices. Then those latter can be replaced
by three new 4-simplices by a move $3\to 3$. After such a
replacement, take in the\/
{\em new\/} matrix $(\partial \omega_i/ \partial L_a)$ the new
submatrix~$\cal B$ containing the same rows and columns, with the only
following change: the row corresponding to the face~$ABC$ is replaced by the
row corresponding to the new face~$DEF$.

The expression
\be
{\displaystyle \det {\cal B}\cdot \prod_{
\hbox{\scriptsize
\begin{tabular}{c}
\rm over all\\[-.5\baselineskip]
\rm 4-simplices
\end{tabular} } } V
\over \displaystyle \prod_{
\hbox{\scriptsize
\begin{tabular}{c}
\rm over all\\[-.5\baselineskip]
\rm 2-dim.\ faces
\end{tabular} } } S }
\label{I.17}
\ee
does not change under such a rebuilding.
}

\medskip

We add here the following. The torsion of complex (\ref{posledT}) is the
alternated product of minors one of which is exactly~$\det {\cal B}$.
With the presented realization of move~$3\to 3$, {\em all other minors\/}
obviously {\em remain unchanged\/}. Thus we can say that formula~(\ref{I.17})
describes the behaviour not only of $\det {\cal B}$ but of the whole
torsion~$\tau$ under a move~$3\to 3$. It is evident also that the acyclicity
property as such, if the complex possessed it, is preserved.

\subsection{Moves $2\leftrightarrow 4$}
\label{subsec 6.2}

For convenience, we speak about moves $2\to 4$, having in mind the obvious
invertibility of the following reasoning and formulas to the case of
moves~$4\to 2$.

Under a move $2\to 4$, one edge and four two-dimensional faces are added
to the simplicial complex. As for the algebraic complex~(\ref{posled}),
the dimensionality of space of column vectors~$(dS_i)$ increases by~4, the
dimensionality of space $(d\vec v_a)$ --- by~3 and the dimensionality of
space~$(d\Omega_a)$ --- by~1. Hence, we can increase by~3 the sizes of the
minor of matrix $(\partial S_i / \partial \vec v_a)$ and by~1 --- those of
the minor of matrix $(\partial \Omega_a / \partial S_i)$, and {\em keep
unchanged the minors\/} of which the torsion is made up {\em beyond the
fragment\/} $(d\Omega_a) \leftarrow (dS_i) \leftarrow (d\vec v_a)$.

As in Section~6 of paper~II, we denote the added edge as $AB$, and the four
added faces will be $ABC$, $ABD$, $ABE$ and~$ABF$. We append the
derivatives of areas of the first three faces to the minor of matrix
$(\partial S_i / \partial \vec v_a)$, while the derivatives
w.r.t.~$S_{ABF}$ --- to the minor of matrix
$(\partial \Omega_a / \partial S_i)$. We denote those minors simply as
$\minor (\partial S_i / \partial \vec v_a)$, etc.: we consider only a single
minor for every matrix, and there is no risk of confusion.

As concerns the 4-simplices, we continue to use the notations of paper~II for
them as well. Namely, under the move $2\to 4$, two adjacent 4-simplices
$ACDEF=\hat B$ and $BCDFE=-\hat A$ are replaced with four ones:
$ABCDE=\hat F$, $ABCFD=-\hat E$, $ABCEF=\hat D$ and~$ABDFE=-\hat C$.

Formula (II.10) shows that under the move $2\to 4$ the important for us minor
of matrix $(\partial S_i / \partial \vec v_a)$ gets multiplied by
\be
\frac{3\, V_{\hat F}\, l_{AB}^5}{S_{ABC}\, S_{ABD}\, S_{ABE}}.
\label{minor-S-v}
\ee
As concerns the minor of matrix $(\partial \Omega_a / \partial S_i)$, it gets
multiplied, according to formula~(II.15), by
\be
\frac{S_{ABF}}{24}\, \frac{V_{\hat A}\, V_{\hat B}}{V_{\hat C}\, V_{\hat D}\, V_{\hat E}}.
\label{minor-L-omega}
\ee
It can be easily derived from here that the following value is conserved
under moves~$2\to 4$:
\be
\frac{\minor\left(\displaystyle \pa{\Omega_a}{S_i} \right)}
{ \minor\left(\displaystyle \pa{S_i}{\vec v_a} \right)}\,
\frac{\displaystyle \prod_{
\hbox{\scriptsize
\begin{tabular}{c}
\rm over all\\[-.5\baselineskip]
\rm 4-simplices
\end{tabular} } } V
\prod_{\hbox{\scriptsize
\begin{tabular}{c}
\rm over all\\[-.5\baselineskip]
\rm edges
\end{tabular} } } 72\,l^5 }
{\displaystyle \prod_{
\hbox{\scriptsize
\begin{tabular}{c}
\rm over all\\[-.5\baselineskip]
\rm 2-dim. faces
\end{tabular} } } S }.
\label{3324}
\ee

Comparing this with formula (\ref{I.17}), we see that the value~(\ref{3324})
is conserved under moves~$3\to 3$ as well (of course,
$\minor (\partial \Omega_a / \partial S_i)$ is the very same thing as
$\det {\cal B}$ in formula~(\ref{I.17})). Besides, it is again evident that
the acyclicity property is conserved.

\section{How the algebraic complex changes under moves $1\leftrightarrow 5$}
\label{sec IIIp7}

Like in the previous Section, we consider, for concreteness, only a move
in one direction, namely $1\to 5$. Under such move, a new {\em vertex\/}
is added to the simplicial complex --- denote it $F$ --- which brings about
the decomposition of one 4-simplex --- denote it $ABCDE$ --- into five
4-simplices, which we denote in the style of the previous Section as
$\hat A$, $\hat B$, $\hat C$, $\hat D$ and~$\hat E$.

We will explain how we extend the minors that enter in the
torsion of complex (\ref{posled}) or~(\ref{posledT}), starting from the left,
that is from the mapping $(dx\hbox{ and }dg)\to (dL_a)$ (the minor of the
``preceding'' mapping $\mathfrak a \to (dx\hbox{ and }dg)$ does not change,
see the remark after formula~(\ref{dL(dxdg)})). The length of columns of
coordinate differentials~$dx$ increases by~4: $dx_F$, $dy_F$, $dz_F$
and~$dt_F$ are added to them --- the differentials of four coordinates of the
new vertex~$F$. Hence, we must choose four~$dL$, corresponding to four new
rows of the minor. Let those be $dL_{AF}$, $dL_{BF}$, $dL_{CF}$ and~$dL_{DF}$.
The following formula holds which can be proved by a direct calculation
using some easy trigonometry (as well as its analogue in 3~dimensions,
see~\cite[formulas (31) and~(32)]{3dim1}):
\be
dx_F \wedge dy_F \wedge dz_F \wedge dt_F = \frac{dL_{AF}\wedge dL_{BF}\wedge
dL_{CF}\wedge dL_{DF}}{384\, V_{\hat E}}\, .
\label{7.1}
\ee
Strictly speaking, here a $\pm$~sign should have been added, but we take
interest in equalities of such kind {\em to within their sign\/}
(as well as in the preceding Section).

It follows from formula~(\ref{7.1}) that the important for us {\em minor
of the matrix of mapping\/} $(dx\hbox{ and }dg)\to (dL_a)$ {\em gets
multiplied by}
\be
384\, V_{\hat E}\, .
\label{dL(dxdg)}
\ee
Note also that the parameters $dg$ responsible for continuous deformations of
representation~$f$ are ``used up'' on algebra~$\mathfrak a$ (if they existed
at all) in the sense that their corresponding rows are included in another
minor, that of mapping $\mathfrak a \to (dx\hbox{ and }dg)$, which does not
change under our move.

Now consider the mapping $(dL_a)\to (d\omega_i)$. Here, only one column is
added to the minor, corresponding to the still ``available''~$dL$,
namely~$dL_{EF}$. On the side of space~$(d\omega_i)$, we add the row
corresponding to~$d\omega_{DEF}$. In consequence, the {\em minor of the matrix
of mapping\/} $(dL_a)\to (d\omega_i)$ {\em gets multiplied by}
\be
\frac{S_{DEF}}{24}\, \frac {V_{\hat E}\, V_{\hat F}}{V_{\hat A}\, V_{\hat B}\, V_{\hat C}}
\label{domegadL}
\ee
(cf.~formula (\ref{minor-L-omega})).

Now we switch from sequence~(\ref{posledT}) to sequence~(\ref{posled}).
Nine rows and columns must be added to the minor of mapping
$(d\vec v_a)\to (dS_i)$. The columns correspond to the nine ``still free''~$dS_i$,
that is $dS_i$ for all added faces~$i$ except $i=DEF$. We are going to
consider them all, step by step, choosing for them on our way nine components
of vectors~$d\vec v_a$ (while their total number is~15: five edges, each
having three deviation components).

Matrix $(\partial S_i / \partial \vec v_a)$ contains many zeroes (the
same applies, by the way, to the other matrices with which we
are occupied here): nonzero entries appear only where
edge~$a$ enters in the boundary of face~$i$. In consequence,
many of its submatrices have a block-triangular form so that their
determinants (i.e., minors of $(\partial S_i / \partial \vec v_a)$) factorize.
In particular, it will be convenient for us to include the derivatives
w.r.t.~{\em all\/} components of deviations of edges $DF$ and~$EF$ in the
minor with which we are occupied now: we will find out that some {\em separate
factors\/} correspond to them as their contributions to the quantity by which
our minor is multiplied. These factors are
\be
\frac{dS_{ADF}\wedge dS_{BDF}\wedge dS_{CDF}}{\left( dv_{DF} \right)_x \wedge
\left( dv_{DF} \right)_y \wedge \left( dv_{DF} \right)_z} =
\frac{3\, V_{\hat E}\, l_{DF}^5}{S_{ADF}\, S_{BDF}\, S_{CDF}}
\label{dSdv1}
\ee
(cf.~formula (\ref{minor-S-v})) for edge~$DF$ and
\be
\frac{dS_{AEF}\wedge dS_{BEF}\wedge dS_{CEF}}{\left( dv_{EF} \right)_x \wedge
\left( dv_{EF} \right)_y \wedge \left( dv_{EF} \right)_z} =
\frac{3\, V_{\hat D}\, l_{EF}^5}{S_{AEF}\, S_{BEF}\, S_{CEF}}
\label{dSdv2}
\ee
for edge~$EF$.

There remain three area differentials --- $dS_{ABF}$, $dS_{ACF}$
and~$dS_{BCF}$ --- for which we must choose three from nine components of
vectors $d\vec v_{AF}$, $d\vec v_{BF}$ and~$d\vec v_{CF}$. To this end, we
first choose three axes for {\em each\/} of these vectors, in order to
take their projections onto these axes as their components. Certainly,
our triples of axes will make up orthonormal coordinate systems in the
three-dimensional spaces orthogonal to corresponding edges. Some easy
reasoning shows that the torsion of complex does not depend on a specific
choice of such coordinate systems.

We require that the $x$~axes (different!) for $d\vec v_{AF}$ and
$d\vec v_{BF}$ lie in the plane~$ABF$ (this fixes their directions, because
they must also be orthogonal to the respective edges). The $y$~axis (common)
for $d\vec v_{AF}$ and~$d\vec v_{BF}$ will lie in the three-dimensional
hyperplane~$ABCF$ (and be, of course, orthogonal to plane~$ABF$). The
$z$~axis --- {\em common for all three\/} $d\vec v_{AF}$, $d\vec v_{BF}$
and~$d\vec v_{CF}$ --- will be orthogonal to hyperplane~$ABCF$.

It remains to choose the directions of axes $x$ and~$y$ for~$d\vec v_{CF}$. It
is enough to fix the direction of the $x$~axis: we choose it to be
orthogonal to plane~$BCF$.

Now we choose three components of vectors $d\vec v$, in order to
include the derivatives with respect to them into the minor of
matrix~$(\partial S_i / \partial \vec v_a)$. These will be $(dv_{BF})_x$,
$(dv_{CF})_x$ and~$(dv_{CF})_y$.

The determinant
\be
\frac{dS_{ABF} \wedge dS_{ACF} \wedge dS_{BCF}}{(dv_{BF})_x \wedge (dv_{CF})_x \wedge (dv_{CF})_y}
\label{3dS3dv}
\ee
again factorizes: it equals the product of the quantity
$\partial S_{ABF} / (\partial v_{BF})_x = l_{BF}$ by the determinant
$\displaystyle \frac{dS_{ACF} \wedge dS_{BCF}}{(dv_{CF})_x \wedge (dv_{CF})_y}$.
But this determinant, too, factorizes due to the fact that $dS_{BCF}$ does not
depend on~$(dv_{CF})_x$ (recall how we chose the direction of axis~$x$
for this vector). Finally we find that expression~(\ref{3dS3dv}) equals
\be
l_{BF}\, l_{CF}^2 \sin \alpha ,
\label{dSdv3}
\ee
where $\alpha$ is the angle between planes $ACF$ and~$BCF$.

Here is the expression --- the product of expressions (\ref{dSdv1}),
(\ref{dSdv2}) and~(\ref{dSdv3}) --- by which the minor of
matrix~$(\partial S_i / \partial \vec v_a)$ is multiplied as a result
of the move $1\to 5$:
\be
\frac{9\, V_{\hat E}\, V_{\hat D}\, l_{DF}^5 \, l_{EF}^5}{S_{ADF}\, S_{BDF}\, S_{CDF}\, S_{AEF}\, S_{BEF}\, S_{CEF}}
\, l_{BF}\, l_{CF}^2 \sin \alpha \, .
\label{dSdv}
\ee

It remains to consider one more minor that changes under the move
$1\to 5$ --- the minor of matrix $(\partial \vec v_a / \partial \sigma)$.
It involves the components $(dv_{AF})_x$, $(dv_{AF})_y$, $(dv_{AF})_z$,
$(dv_{BF})_y$, $(dv_{BF})_z$ and~$(dv_{CF})_z$ of edge deviations and,
on the other hand, all six components of bivector~$d\sigma$ corresponding
to vertex~$F$.

Bivector~$d\sigma$ can be thought of as an element of Lie algebra
$\mathfrak s \mathfrak o (4)$, and its components --- as infinitesimal
rotation angles within the six coordinate planes. We have to choose a proper
coordinate system for it. This time, we will denote its axes by numbers
1, 2, 3, 4 (rather than letters $x$, $y$, $z$,~$t$).

We choose axes 1 and 2 to lie in the plane $ABF$: axis~1 will be orthogonal
to vector~$\overrightarrow {BF}$, and axis~2 --- parallel to it. Axis~3
will lie in the hyperplane~$ABCF$ (and, of course, orthogonal to $ABF$), and
axis~4 --- orthogonal to~$ABCF$.

Now (the thing we are already accustomed to) the value by which the
minor is multiplied again factorizes. First (cf.~formula~(\ref{sigma->v})),
the components $(dv_{AF})_z$, $(dv_{BF})_z$ and~$(dv_{CF})_z$ are determined
just by three rotation angles in the direction of axis~4, i.e.\
$d\sigma_{14}$, $d\sigma_{24}$ and~$d\sigma_{34}$. Namely, they are connected
by means of the following $3\times 3$ submatrix of matrix
$(\partial \vec v_a / \partial \sigma)$ (see again~(\ref{sigma->v})):
\be
\left( \begin{array}{ccc}
\frac{(AF)_1}{L_{AF}} & \frac{(AF)_2}{L_{AF}} & \frac{(AF)_3}{L_{AF}} \\
\noalign{\medskip}
\frac{(BF)_1}{L_{BF}} & \frac{(BF)_2}{L_{BF}} & \frac{(BF)_3}{L_{BF}} \\
\noalign{\medskip}
\frac{(CF)_1}{L_{CF}} & \frac{(CF)_2}{L_{CF}} & \frac{(CF)_3}{L_{CF}}
\end{array} \right) ,
\label{4z}
\ee
where, for instance, $(AF)_1$ is the component of
vector~$\overrightarrow{AF}$ along axis~1. The determinant of
matrix~(\ref{4z}) is
\be
\frac{6\, V_{ABCF}}{L_{AF}\, L_{BF}\, L_{CF}} =
\frac{4\, S_{ACF}\, S_{BCF} \sin \alpha}{L_{AF}\, L_{BF}\, l_{CF}^3},
\label{dvdsigma1}
\ee
where $V_{ABCF}$ is the {\em three-dimensional\/} volume, while $\alpha$ is
the same angle as in formula~(\ref{dSdv3}).

There remain three rotations within hyperplane $ABCF$, and three components
$(dv_{AF})_x$, $(dv_{AF})_y$ and~$(dv_{BF})_y$. Here, too, the factorizability
persists: the rotation within the plane~$ABF$ affects only~$(dv_{AF})_x$, and
this gives the factor
\be
\frac{(\partial v_{AF})_x}{\partial \sigma_{12}} = \frac {1}{l_{AF}}\, .
\label{dvdsigma2}
\ee
The rotation within the plane orthogonal to~$BF$ affects only~$(dv_{AF})_y$,
which gives the factor
\be
\frac{(\partial v_{AF})_y}{\partial \sigma_{13}} = \frac {\sin \beta}{l_{AF}}\, ,
\label{dvdsigma3}
\ee
$\beta$ being the angle between edges $AF$ and~$BF$. Now the last rotation
remains with the corresponding factor
\be
\frac{(\partial v_{BF})_y}{\partial \sigma_{23}} = \frac {1}{l_{BF}}\, .
\label{dvdsigma4}
\ee

Taking the product of the right-hand sides of formulas (\ref{dvdsigma1}),
(\ref{dvdsigma2}), (\ref{dvdsigma3}) and~(\ref{dvdsigma4}),
we find the factor by which the minor of matrix
$(\partial \vec v_a / \partial \sigma)$ gets multiplied under our
move~$1\to 5$:
\be
\frac{8\, S_{ACF}\, S_{BCF}\, S_{ABF}}{l_{AF}^5\, l_{BF}^4\, l_{CF}^3}
\cdot \sin \alpha\, ,
\label{dvdsigma}
\ee
where we have taken into account that
$\sin \beta = 2\, S_{ABF} / l_{AF}\, l_{BF}$.

In the same way as in Subsections \ref{subsec 6.1} and~\ref{subsec 6.2},
it is clear from the performed reasoning and calculations that the property
of acyclicity of the complex is conserved under the considered moves.

\section{The formula for the manifold invariant}
\label{sec IIIp8}

In the two preceding Sections we have analyzed the behaviour of the minors
whose alternated product makes up the torsion of the complex (\ref{posled})
or~(\ref{posledT}) (if this complex is acyclic), under Pachner moves
$3\to 3$, $2\leftrightarrow 4$ and~$1\leftrightarrow 5$. We have now to unite
these results to obtain the quantity which {\em does not depend on a
triangulation}, i.e.\ an {\em invariant of a four-dimensional
piecewise-linear manifold}.

We choose the signs in the alternated product in such way that our formula for
the invariant look as similar as possible to the ``three-dimensional''
formula~(II.5). Namely, we take the minor of the matrix of mapping
$\mathfrak a \to (dx\hbox{ and }dg)$ raised in the power~$-1$, then the minor
of the matrix of mapping $(dx\hbox{ and }dg) \to (dL_a)$ raised in the
power~$+1$ and so~on. We will get the factor by which the so defined
torsion~$\tau$ is multiplied under the move~$1\to 5$ when we multiply the
expressions (\ref{dL(dxdg)}) and~(\ref{dSdv}) and divide by the product of
(\ref{domegadL}) and~(\ref{dvdsigma}). The result can be written as
$$
2^7 \cdot 3^4 \cdot \frac{V_{\hat A}\, V_{\hat B}\, V_{\hat C}\, V_{\hat D}\, V_{\hat E}}{V_{\hat F}}
\cdot \frac{\displaystyle \prod_{
\hbox{\scriptsize
\begin{tabular}{c}
\rm over new\\[-.5\baselineskip]
\rm edges
\end{tabular} } } l^5 }{\displaystyle \prod_{
\hbox{\scriptsize
\begin{tabular}{c}
\rm over new\\[-.5\baselineskip]
\rm 2-dim.~faces
\end{tabular} } } S } \; .
$$

Comparing this with the results of Section~\ref{sec IIIp6} (formulas
(\ref{I.17}) and~(\ref{3324})), we find the following final expression for the
invariant of a four-dimensional manifold in terms of the torsion of complex
(\ref{posled}) or~(\ref{posledT}) and other Euclidean geometric values:
\be
I=2^{-16} \cdot 3^{-12} \cdot
\frac{\displaystyle \tau \cdot \prod_{
\hbox{\scriptsize
\begin{tabular}{c}
\rm over all\\[-.5\baselineskip]
\rm 2-dim.~faces
\end{tabular} } } S }{\displaystyle \prod_{
\hbox{\scriptsize
\begin{tabular}{c}
\rm over all\\[-.5\baselineskip]
\rm 4-simplices
\end{tabular} } } V \prod_{
\hbox{\scriptsize
\begin{tabular}{c}
\rm over all\\[-.5\baselineskip]
\rm edges
\end{tabular} } } 72\, l^5
} \cdot \left( 2^8 \cdot 3^6 \right)^{(\hbox{\scriptsize number of vertices})} \, .
\label{i4}
\ee
The factor $2^{-16} \cdot 3^{-12}$ has been added in order that the invariant
be equal to unity for the sphere~$S^4$, see the following Section.

\section{Example: sphere~$S^4$}
\label{sec IIIp9}

The fundamental group of sphere~$S^4$ is trivial, thus algebra~$\mathfrak a$
in sequence~(\ref{posledT}) will be the whole Lie algebra~$\mathfrak e_4$ of
motions of the four-dimensional Euclidean space. There are obviously no
continuous parameters~$dg$ describing the deformations of representation
$\pi_1(S^4) \to E_4$. The space of ``trivial'' vertex deviations is
ten-dimensional and consists of deviations of type~(\ref{trivsigma}) with
arbitrary $d\tau_{\alpha \beta}$ and~$ds_{\beta}$.

We take the canonical triangulation of~$S^4$ consisting of two 4-simplices,
both with vertices $A$, $B$, $C$, $D$ and~$E$. For such triangulation, all
deficit angles $d\omega_i$ and~$d\Omega_a$ are identically equal to zero,
because any~$\omega_i$ is obtained by summing two terms differing only in sign
(the $+$~sign is ascribed to one of 4-simplices, while the $-$~sign to the other one,
see~I, Section~3). Thus, matrix $(\partial \Omega_a / \partial S_i)$ and its
conjugate~$(\partial \omega_i / \partial L_a)$ are zero. This means that the
torsion of complex (\ref{posled}) or~(\ref{posledT}) factorizes in the
product over two following sequences:
\be
0 \rightarrow \mathfrak e_4 \rightarrow (dx) \rightarrow (dL_a) \rightarrow 0
\label{l}
\ee
and
\be
0 \leftarrow (dS_i) \leftarrow (d\vec v_a) \leftarrow (d\sigma) \leftarrow (d\sigma)_0 \leftarrow 0.
\label{pr}
\ee

We fix a Euclidean coordinate system in the space~$\mathbb R_4$, with axes
$x,y,z,t$. We place the vertex~$A$ of our triangulation into the origin of
coordinates, and the remaining vertices --- in points $B(1,0,0,0)$,
$C(0,1,0,0)$, $D(0,0,1,0)$ and~$E(0,0,0,1)$.

\subsection{Calculations for sequence (\protect\ref{l})}

The ten-dimensional algebra $\mathfrak e_4$ consists of infinitesimal
rotations and translations (of course, we measure the rotations within the six
coordinate planes in radians, and the translations --- in the units of
coordinate axes), whereas the twenty-dimensional space~$(dx)$ --- of column
vectors $\pmatrix{dx_A \cr \vdots \cr dt_E}$. We take the minor of matrix
of mapping $\mathfrak e_4 \to (dx)$ corresponding to the following
{\em ten\/} coordinate differentials: $dx_A, dy_A, dz_A,\allowbreak
dt_A, dy_B, dz_B,\allowbreak dt_B, dz_C, dt_C, dt_D.$

 From this minor, a unity factor splits off at once corresponding to the
mapping (translations) $\to$ (vertex~$A$ coordinates); this splitting off is
caused by the fact that {\em rotations do not affect\/} the coordinates
of~$A$. There remains the mapping of six rotations into sets
$(dy_B,dz_B,dt_B,\allowbreak dz_C,dt_C,dt_D)$. Factorization still works
here: rotations within the planes $xy$, $xz$ and~$xt$ affect only point~$B$
(from those remaining) and so on. The result is: {\em the minor of mapping
$\mathfrak e_4 \to (dx)$ equals unity\/} (as usual, to within a sign).

Now we consider the minor of matrix of the mapping $(dx)\to (dL_a)$
corresponding to the set $(dx_B, dx_C, dy_C,\allowbreak
dx_D, dy_D, dz_D,\allowbreak dx_E, dy_E, dz_E, dt_E)$.
It, too, factorizes into the product of determinants of mappings
\begin{eqnarray}
(dx_E, dy_E, dz_E, dt_E) & \to & (dL_{AE}, dL_{BE}, dL_{CE}, dL_{DE}),
\label{9.1.1} \\
(dx_D, dy_D, dz_D) & \to & (dL_{AD}, dL_{BD}, dL_{CD}),
\label{9.1.2} \\
(dx_C, dy_C) & \to & (dL_{AC}, dL_{BC}),
\label{9.1.3} \\
dx_B & \to & dL_{AB}.
\label{9.1.4}
\end{eqnarray}
Recall that $L$ is the squared length of a corresponding edge. A direct
calculation shows that the determinants of mappings
(\ref{9.1.1}), (\ref{9.1.2}), (\ref{9.1.3}) and~(\ref{9.1.4}) equal
$2^4$, $2^3$, $2^2$ and~$2$, respectively.

{\bf Conclusion: }the multiplicative contribution of sequence (\ref{l}) to
the torsion is~$2^{10}$.

\subsection{Calculations for sequence (\protect\ref{pr})}

We start from the right, i.e.\ from the mapping $(d\sigma)_0 \to (d\sigma)$.
The ten-dimensional space~$(d\sigma)_0$ consists of the components of
antisymmetric tensor~$d\tau_{\alpha \beta}$ and vector~$ds_{\beta}$, whereas
the thirty-dimensional space~$(d\sigma)$ --- of the components of all five
vertices' deviations. The mapping is given by formula~(\ref{trivsigma}).

For composing the minor of matrix of the mapping $(d\sigma)_0 \to (d\sigma)$, we choose
the following ten components of tensors $d\sigma_A,\ldots,d\sigma_E$: we take
{\em all six\/} components of~$d\sigma_A$ and {\em one\/} component for each
of the remaining deviations, namely, $(d\sigma_B)_{zt}$, $(d\sigma_C)_{zt}$,
$(d\sigma_D)_{xy}$ and~$(d\sigma_E)_{xy}$. It makes sense to remind here that
this choice depends on us (the only requirement is that minors be nonzero),
and this exactly choice was motivated by the convenience of further
calculations.

The vector~$ds_{\beta}$ does not affect $d\sigma_A$, because the radius
vector of point~$A$ is zero. Hence, our minor factorizes, and the factor
corresponding to mapping $d\tau \to d\sigma_A$ is unity (because this is an
identical mapping: $d\sigma_A=d\tau$). The remaining $4\times 4$ minor
factorizes into the product of four unit factors, because, for instance,
$(d\sigma_B)_{zt}$ depends on~$ds_y$ only (this can be seen from
formula~(\ref{trivsigma}) if we replace in it the subscript $A$ with~$B$,
ignore $d\tau_{\alpha \beta}$ and observe that the single nonzero component
of vector~$\vec r_B$ is $(r_B)_x$) and so~on.

Thus, our selected minor of matrix of the mapping $(d\sigma)_0 \to (d\sigma)$
turned out to equal unity.

We pass on to the mapping $(d\sigma)\to (d\vec v_a)$. Here one should fix at
first a three-dimensional basis for each of ten vectors
$d\vec v_{AB},\ldots,d\vec v_{DE}$. It is done most easily for those edges
that begin at point~$A$: edge~$AB$ lies on the $x$~axis, so we choose as the
coordinate axes for $d\vec v_{AB}$ the three remaining axes $y$, $z$ and~$t$;
then we choose axes in a similar way for the deviations of edges $AC$, $AD$
and~$AE$.

Each of the remaining six edges lies within some coordinate {\em plane}, and
we will treat them in the following way. Edge~$BC$ lies in the plane~$xy$; we
choose the following three axes for~$d\vec v_{BC}$: the bisector of the angle
formed by axes $x$ and~$y$, and also axes $z$ and~$t$; we follow this model
in choosing the axes for deviations of edges $BD$, $BE$, $CD$, $CE$ and~$DE$.
We will denote the components of vectors of edge deviations along the
bisectors of coordinate angles as $\left( dv_{BC} \right)_{\rm bis}$, etc.
(keeping in mind that every deviation has its own bisector for an axis).

We must choose twenty components of vectors of edge deviations (because we
have twenty components of vertex deviations not included in our previous
minor). We choose {\em all\/} components for edges $AB$, $AC$, $AD$ and~$AE$,
and also the following eight components: $(dv_{BD})_y$, $(dv_{BD})_t$,
$(dv_{BE})_y$, $(dv_{BE})_z$, $(dv_{CD})_x$, $(dv_{CD})_t$, $(dv_{CE})_x$
and~$(dv_{CE})_z$. We drew in Figure~\ref{IIIris2}
\begin{figure}
\begin{center}
\unitlength=1.00mm
\special{em:linewidth 0.6pt}
\linethickness{0.6pt}
\begin{picture}(97.00,67.00)
\emline{5.00}{35.00}{1}{50.00}{5.00}{2}
\emline{50.00}{5.00}{3}{95.00}{35.00}{4}
\emline{95.00}{35.00}{5}{50.00}{65.00}{6}
\emline{50.00}{65.00}{7}{5.00}{35.00}{8}
\put(3.00,35.00){\makebox(0,0)[rc]{$B$}}
\put(50.00,3.00){\makebox(0,0)[ct]{$E$}}
\put(97.00,35.00){\makebox(0,0)[lc]{$C$}}
\put(50.00,67.00){\makebox(0,0)[cb]{$D$}}
\put(50.00,63.50){\makebox(0,0)[ct]{$\scriptstyle z$}}
\put(92.50,35.00){\makebox(0,0)[rc]{$\scriptstyle y$}}
\put(50.00,6.50){\makebox(0,0)[cb]{$\scriptstyle t$}}
\put(7.50,35.00){\makebox(0,0)[lc]{$\scriptstyle x$}}
\put(5.00,29.00){\makebox(0,0)[cc]{$d\sigma_{yt}$}}
\put(5.00,41.00){\makebox(0,0)[cc]{$d\sigma_{yz}$}}
\put(41.00,5.00){\makebox(0,0)[cc]{$d\sigma_{zx}$}}
\put(41.00,65.00){\makebox(0,0)[cc]{$d\sigma_{tx}$}}
\put(26.00,15.00){\makebox(0,0)[cc]{$dv_z$}}
\put(20.00,19.00){\makebox(0,0)[cc]{$dv_y$}}
\put(59.00,5.00){\makebox(0,0)[cc]{$d\sigma_{zy}$}}
\put(95.00,29.00){\makebox(0,0)[cc]{$d\sigma_{xt}$}}
\put(95.00,41.00){\makebox(0,0)[cc]{$d\sigma_{xz}$}}
\put(59.00,65.00){\makebox(0,0)[cc]{$d\sigma_{ty}$}}
\put(20.00,51.00){\makebox(0,0)[cc]{$dv_y$}}
\put(26.00,55.00){\makebox(0,0)[cc]{$dv_t$}}
\put(74.00,55.00){\makebox(0,0)[cc]{$dv_t$}}
\put(80.00,51.00){\makebox(0,0)[cc]{$dv_x$}}
\put(80.00,19.00){\makebox(0,0)[cc]{$dv_x$}}
\put(74.00,15.00){\makebox(0,0)[cc]{$dv_z$}}
\special{em:linewidth 0.4pt}
\linethickness{0.4pt}
\put(36.50,8.00){\vector(-3,2){6.00}}
\put(9.50,26.00){\vector(3,-2){6.00}}
\put(9.50,44.00){\vector(3,2){6.00}}
\put(36.50,62.00){\vector(-3,-2){6.00}}
\put(63.50,62.00){\vector(3,-2){6.00}}
\put(90.50,44.00){\vector(-3,2){6.00}}
\put(90.50,26.00){\vector(-3,-2){6.00}}
\put(63.50,8.00){\vector(3,2){6.00}}
\put(50.00,5.00){\circle*{1.00}}
\put(95.00,35.00){\circle*{1.00}}
\put(50.00,65.00){\circle*{1.00}}
\put(5.00,35.00){\circle*{1.00}}
\end{picture}
\end{center}
\caption{}
\label{IIIris2}
\end{figure}
the four edges to which these components belong, and wrote out these
components near the edges (for example, we wrote $dv_y$ and~$dv_t$ near the
edge~$BD$ meaning $(dv_{BD})_y$ and~$(dv_{BD})_t$).

Of course, our $20\times 20$ minor again greatly factorizes. For example, the minor
corresponding, on the one hand, to the components $(d\sigma_B)_{xy}$,
$(d\sigma_B)_{xz}$ and~$(d\sigma_B)_{xt}$, and on the other hand --- to the
three components of~$d\vec v_{AB}$, factors out (this is because all the
components of~$d\sigma_A$ are already used up, whereas only they of all the
remaining components of~$d\sigma$ could influence~$d\vec v_{AB}$). From
formulas of type~(\ref{sigma->v}) we can see that this factor equals unity.
Similarly, three more minors factor out which are obtained from the above minor by
changing $B\to C$, $x\leftrightarrow y$, or by changing $B\to D$,
$x\leftrightarrow z$, or by changing $B\to E$, $x\leftrightarrow t$.

After this, there remain the eight already mentioned components of
vectors~$d\vec v$ depicted in Figure~\ref{IIIris2}, and the eight
components of~$d\sigma$ also depicted in Figure~\ref{IIIris2}, with the
understanding that if, e.g., $d\sigma_{yz}$ and~$d\sigma_{yt}$ are drawn near
the vertex~$B$, then they are $(d\sigma_B)_{yz}$ and~$(d\sigma_B)_{yt}$.
The corresponding minor factorizes in eight separate factors.
The reason for
this is that each component of $d\vec v$ depends on {\em only one\/} component
of~$d\sigma$. These dependencies are shown by arrows in Figure~\ref{IIIris2}.
Note also that the small letters near the vertices denote the axes where these
vertices lie.

We see from formulas of type (\ref{sigma->v}) or~(\ref{4.1}) that all the
eight factors equal~$\frac{1}{2}$. Consequently, their contribution to the 
torsion is~$2^8$.

It remains for us to consider the mapping $(d\vec v_a)\to (dS_i)$. As we
remember, all the components of vectors~$d\vec v$ corresponding to edges
$AB$, $AC$, $AD$ and~$AE$ are already used up. For the remaining six edges,
we have the components of~$d\vec v$ along the bysectors of coordinate angles,
as well as $(dv_{BC})_z$, $(dv_{BC})_t$, $(dv_{DE})_x$ and~$(dv_{DE})_y$
(these latter belong to the two edges {\em absent\/} from Figure~\ref{IIIris2}).

Each of the six differentials~$dS_i$, where face~$i$ contains vertex~$A$, is
influenced by only one component of only one of the remaining 
vectors~$d\vec v$. In this way, six factors appear: 
$\partial S_{ABC} / (\partial v_{BC})_{\rm bis} =\, \root\of{2}$ and five
other ones, all equal to it. So, here we have the contribution to the
torsion equal to~$2^3$.

There remains the minor whose rows correspond to componentÁs
$dS_{BCD}$, $dS_{BCE}$, $dS_{BDE}$ and~$dS_{CDE}$, while columns --- to
$(dv_{BC})_z$, $(dv_{BC})_t$, $(dv_{DE})_x$ and~$(dv_{DE})_y$. Here, too,
each component~$dS$ turns out to depend on only one component~$dv$.
We get the product of four partial derivatives, all equal to each other;
for example, one of them is
$$ 
\frac{\partial S_{BCD}}{(\partial
v_{BC})_z} = l_{BC} \cdot \cos \gamma =\, \root\of{2}
 \cdot \,
\root\of{\frac{2}{3}} = \frac{2}{\root\of{3}}, 
$$ 
where $\gamma$ is the angle between face $BCD$ and axis~$z$. Thus, here the
contribution to the torsion is~$2^4\cdot 3^{-2}$.

{\bf Conclusion: }sequence (\ref{pr}) makes the multiplicative contribution
to the torsion, equal to $2^8 \cdot 2^3 \cdot 2^4
\cdot 3^{-2} = 2^{15} \cdot 3^{-2}$.

\subsection{The result: invariant for sphere~$S^4$}

Combining the conclusions made in the end of two previous subsections, we find
that the torsion for sphere~$S^4$ is
$$
\tau (S^4) = 2^{25}\cdot 3^{-2}.
$$
Now we calculate the products entering in formula~(\ref{i4}).

In our complex, there are six two-dimensional faces of area~$1/2$ and four
faces of area~$\,\root\of{3}/2$. Next, $\prod S = 3^2 / 2^{10}$.
Then, there are two 4-simplices, both of volume~$1/24$, thus 
$\prod V = 2^{-6}\cdot 3^{-2}$. Finally, there are four edges of length~$1$ and
six edges of length~$\,\root\of{2}$, thus 
$\prod 72 \, l^5 = 2^{45}\cdot 3^{20}$.

All this together leads to the formula announced in the end of 
Section~\ref{sec IIIp8}:
\be
I(S^4) = 1.
\ee

\section{Discussion}
\label{sec IIIobsuzhdenie}

So, in the case of sphere~$S^4$ the complex turned out to be acyclic, and we
managed to calculate its torsion (and our invariant). The largest
determinant that we had to deal with was of sizes $20\times 20$, but, luckily,
it factorized in a product of smaller determinants.

Hopefully, new properties of our invariants will be discovered with time,
which will simplify the calculations, and some relevant techniques will be
elaborated. This will give the real possibility to calculate the invariants
for a large enough manifold zoo. At this moment, one of the interesting
questions is what we will get for the product of two-dimensional 
spheres~$S^2\times S^2$ and whether we will be able to do something if the
corresponding complex turns out {\em not\/} to be acyclic.

One more problem is the generalization of our complexes and finding their
possible quantum analogues. One can begin with constructing a complex based on
the $SL(2)$-solution to the pentagon equation from paper~\cite{KM2}.

It looks quite plausible that our constructions can be generalized in such way that
they include also the {\em Reidemeister torsion}. In prospect, one can think
about the creation of a new general theory that will combine the ideas of the
algebra of acyclic complexes and quantum topology.

\medskip

{\bf Acknowledgements. }\vadjust{\nobreak}The work has been performed with the
partial financial support from Russian Foundation for Basic Research, Grant
no.~01-01-00059.

\end{document}